\documentclass{amsart}
\usepackage[utf8]{inputenc}
\usepackage{amsmath,amsthm}
\usepackage{tensor}
\usepackage{amssymb}
\usepackage{stmaryrd}
\usepackage{color}
\usepackage[left=3cm,right=3cm,top=3cm,bottom=4.5cm]{geometry}
\usepackage{pdfsync}
\usepackage{wasysym}
\usepackage{tikz}
\usetikzlibrary{decorations.pathmorphing,cd,patterns}
\usepackage{url}
\usepackage{enumitem}
\usepackage{hyperref}
\hypersetup{linktoc=all}
\usepackage[capitalise,noabbrev]{cleveref}

\theoremstyle{plain}
	\newtheorem{thm}{Theorem}[section]
	\newtheorem*{thm*}{Theorem}
	\newtheorem{cor}[thm]{Corollary}
	\newtheorem*{cor*}{Corollary}

	\newtheorem*{prop*}{Proposition}
	\newtheorem{lem}[thm]{Lemma}
	\newtheorem{lemma}[thm]{Lemma}
	\newtheorem*{lem*}{Lemma}
	
	\newtheorem*{ex*}{Exercise}
	
	\newtheorem*{claim*}{Claim}
	
	\newtheorem*{question*}{Question}
	
	\newtheorem*{fact*}{Fact}
\theoremstyle{definition}
	
	\newtheorem{definition}[thm]{Definition}
	\newtheorem*{Def*}{Definition}
	
	\newtheorem*{obs*}{Observation}
	\newtheorem{rmk}[thm]{Remark}
	\newtheorem*{rmk*}{Remark}
	
	\newtheorem{soln*}{Solution}
	
	\newtheorem*{note*}{Note}
	
	\newtheorem{example}[thm]{Example}
	\newtheorem*{eg*}{Example}	
	
	\newtheorem*{construction*}{Construction}
	
	\newtheorem*{warning*}{Warning}
	
	\newtheorem*{conj*}{Conjecture}
	
	\newtheorem*{recall*}{Recall}
	
	\newtheorem*{convention*}{Convention}
	
\newcommand{\nats}{\mathbb{N}}

\newcommand{\op}{\mathrm{op}}
\newcommand{\id}{\mathrm{id}}

\newcommand{\Hom}{\mathrm{Hom}}

\newcommand{\Ob}{\operatorname{Ob}}

\newcommand{\Psh}{\mathsf{Psh}}
\newcommand{\Fun}{\mathrm{Fun}}
\DeclareMathOperator{\Nat}{Nat}

\newcommand{\Set}{\mathsf{Set}}

\newcommand{\Cat}{\mathsf{Cat}}

\newcommand{\calA}{\mathcal{A}}
\newcommand{\calB}{\mathcal{B}}
\newcommand{\calC}{\mathcal{C}}
\newcommand{\calD}{\mathcal{D}}
\newcommand{\calE}{\mathcal{E}}

\newcommand{\calK}{\mathcal{K}}

\newcommand{\calO}{\mathcal{O}}

\newcommand{\pre}{\preceq}
\newcommand{\eepsilon}{{\boldsymbol\epsilon}}
\newcommand{\zzeta}{{\boldsymbol\zeta}}
\newcommand{\eeta}{{\boldsymbol\eta}}
\newcommand{\ddelta}{{\boldsymbol\delta}}
\usepackage[capitalise]{cleveref}
\newcommand{\nq}{{[n;\mathbf{q}]}}
\mathchardef\mhyphen="2D
\newcommand{\twoCat}{2\mhyphen\Cat}
\newcommand{\inftytwoCat}{(\infty,2)\mhyphen\Cat}
\newcommand{\inftynCat}{(\infty,n)\mhyphen\Cat}
\newcommand{\Thetahat}{\widehat{\Theta_2}}
\newcommand{\nerve}{\operatorname{N}_\Theta}

\DeclareFontFamily{U}{min}{}
\DeclareFontShape{U}{min}{m}{n}{<-> udmj30}{}
\newcommand\yo{\!\text{\usefont{U}{min}{m}{n}\symbol{'207}}\!}

\newcommand{\Spaces}{\mathsf{Spaces}}
\newcommand{\nerveinfty}{\operatorname{N}_\square}
\newcommand{\Alg}{\operatorname{Alg}}
\newcommand{\strong}{\mathrm{strong}}
\newcommand{\lax}{\mathrm{lax}}
\newcommand{\cocts}{\mathrm{cocts}}
\newcommand{\refl}{L}

\title{A model-independent Gray tensor product for $(\infty,2)$-categories}
\author{Timothy Campion and Yuki Maehara}
\date{\today}

\begin{document}

\maketitle

\begin{abstract}
    We construct a (lax) Gray tensor product of $(\infty,2)$-categories and characterize it via a model-independent universal property. Namely, it is the unique monoidal biclosed structure on the $\infty$-category of $(\infty,2)$-categories which agrees with the classical Gray tensor product of strict 2-categories when restricted to the Gray cubes (i.e. the Gray tensor powers $[1]^{\otimes n}$ of the arrow category).
\end{abstract}

\tableofcontents

\section{Introduction}
The Gray tensor product was first introduced by Gray \cite{Gray} in the setting of of strict 2-categories. It has since been extended to strict $\omega$-categories as well as weak $(\infty,\infty)$-categories and in $(\infty,2)$-categories by various authors in various models (we review the literature in detail shortly). Although there exist various comparison results among these models, the proliferation of ``Gray tensor products" begs the question of how to characterize the Gray tensor product more abstractly. In this note, we address this question by showing (\cref{thm:day-refl}) that the Gray tensor product of the second author \cite{Maehara:Gray} in fact descends from the 2-quasicategory model to a monoidal structure on the $\infty$-category $\inftytwoCat$ of $(\infty,2)$-categories. Moreover, we show (\cref{cor:univ-prop}) that this ``Gray tensor product" enjoys a universal property: it is the unique monoidal biclosed structure on $\inftytwoCat$ which agrees with the classical Gray tensor product on the Gray cubes (i.e. the Gray tensor powers $\square^n = [1]^{\otimes n}$ of the arrow category). In this way, we reduce the specification of the Gray tensor product on $\inftytwoCat$ to the specification of a small amount of combinatorial data.

\begin{figure}
\[
\begin{tikzpicture}[scale = 1.5, baseline = -2]
	\filldraw
	(0,0) circle [radius = 1pt]
	(1,0) circle [radius = 1pt];
	
	\draw[->] (0.1,0) -- (0.9,0);

\end{tikzpicture} \hspace{50pt}
\begin{tikzpicture}[scale = 1.5, baseline = 18.5]
	\filldraw
	(0,0) circle [radius = 1pt]
	(1,0) circle [radius = 1pt]
	(0,1) circle [radius = 1pt]
	(1,1) circle [radius = 1pt];
	
	\draw[->] (0.1,0) -- (0.9,0);
	\draw[->] (0.1,1) -- (0.9,1);
	\draw[<-] (0,0.1) -- (0,0.9);
	\draw[<-] (1,0.1) -- (1,0.9);
	
	\draw[->, double] (0.7,0.7) -- (0.3,0.3);
\end{tikzpicture} \hspace{50pt}
\begin{tikzpicture}[baseline = -2,scale = 1.5]	
	\filldraw
	(150:1) circle [radius = 1pt]
	(90:1) circle [radius = 1pt]
	(30:1) circle [radius = 1pt]
	(-30:1) circle [radius = 1pt]
	(-90:1) circle [radius = 1pt]
	(-150:1) circle [radius = 1pt]
	(0,0) circle [radius = 1pt];
	
	\draw[->] (150:1) + (30:0.1) --+ (30:0.9);
	\draw[->] (90:1) + (-30:0.1) --+ (-30:0.9);
	\draw[->] (30:1) + (-90:0.1) --+ (-90:0.9);
	\draw[->] (150:1) + (-90:0.1) --+ (-90:0.9);
	\draw[->] (-150:1) + (-30:0.1) --+ (-30:0.9);
	\draw[->] (-90:1) + (30:0.1) --+ (30:0.9);

	\draw[->] (150:1) + (-30:0.1) --+ (-30:0.9);
	\draw[->] (0:0) + (-90:0.1) --+ (-90:0.9);
	\draw[->] (0:0) + (30:0.1) --+ (30:0.9);
	
	\draw[<-,double] (-150:0.5) + (-0.15,-0.15) --+ (0.15,0.15);
	\draw[<-,double] (90:0.5) + (-0.15,-0.15) --+ (0.15,0.15);
	\draw[<-,double] (-30:0.5) + (-0.15,-0.15) --+ (0.15,0.15);
\end{tikzpicture}
\quad = \quad
\begin{tikzpicture}[baseline = -2,scale = 1.5]
	\filldraw
	(150:1) circle [radius = 1pt]
	(90:1) circle [radius = 1pt]
	(30:1) circle [radius = 1pt]
	(-30:1) circle [radius = 1pt]
	(-90:1) circle [radius = 1pt]
	(-150:1) circle [radius = 1pt]
	(0,0) circle [radius = 1pt];
	
	\draw[->] (150:1) + (30:0.1) --+ (30:0.9);
	\draw[->] (90:1) + (-30:0.1) --+ (-30:0.9);
	\draw[->] (30:1) + (-90:0.1) --+ (-90:0.9);
	\draw[->] (150:1) + (-90:0.1) --+ (-90:0.9);
	\draw[->] (-150:1) + (-30:0.1) --+ (-30:0.9);
	\draw[->] (-90:1) + (30:0.1) --+ (30:0.9);

	\draw[->] (0:0) + (-30:0.1) --+ (-30:0.9);
	\draw[->] (90:1) + (-90:0.1) --+ (-90:0.9);
	\draw[->] (-150:1) + (30:0.1) --+ (30:0.9);
	
	\draw[<-,double] (150:0.5) + (-0.15,-0.15) --+ (0.15,0.15);
	\draw[<-,double] (-90:0.5) + (-0.15,-0.15) --+ (0.15,0.15);
	\draw[<-,double] (30:0.5) + (-0.15,-0.15) --+ (0.15,0.15);
\end{tikzpicture}
\]
\caption{$\square^1 = [1]$, $\square^2 = [1] \otimes [1]$, and $\square^3 = [1] \otimes [1] \otimes [1]$}\label{cubes-0}
\end{figure}
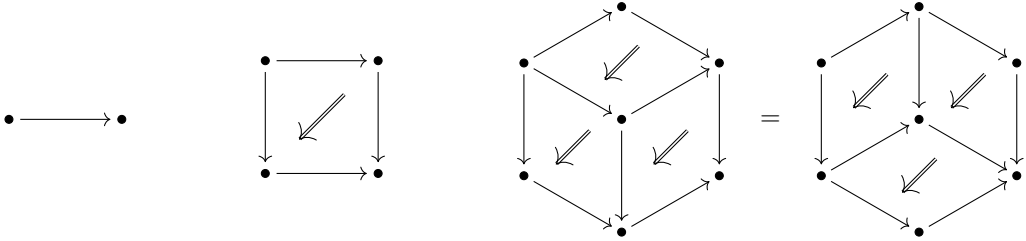

Let us put this result into context. The Gray tensor product of strict 2-categories was introduced in \cite{Gray} in order to study lax natural transformations. It is a monoidal biclosed structure $\otimes$ on the 1-category $\twoCat$ of strict 2-categories, such that the left (resp. right) internal hom $\llbracket A, B\rrbracket$ has strict 2-functors $A \to B$ for objects, lax (resp. oplax) natural transformations for 1-morphisms, and modifications for 2-morphisms. The unit of the monoidal structure is the terminal category $\square^0$, and $\Ob(A \otimes B) = \Ob A \times \Ob B$. The first few Gray tensor powers of the arrow category $[1]$ are pictured in \cref{cubes-0}. In general, the power $[1]^{\otimes n}$ is a lax $n$-dimensional cube. Gray's construction was extended to strict $\omega$-categories using pasting schemes \cite{Crans:thesis}, cubical $\omega$-categories \cite{al-Agl;Brown;Steiner:Multiple} and augmented directed complexes \cite{Steiner:omega}.

Several authors have constructed related monoidal structures on (weak) $(\infty,n)$-categories for general $n \in \nats \cup \{\infty\}$. Verity's complicial model \cite{Verity:I} admits a monoidal structure and a Quillen bifunctor which are weakly equivalent and homotopically well-behaved. The two versions of tensor product agree when reflected to a suitable subcategory, giving rise to a monoidal biclosed structure on $\inftynCat$ (see \cite[Theorem 1.33]{Campion;Kapulkin;Maehara}). The comical models (\cite{Campion;Kapulkin;Maehara}, \cite{Doherty;Kapulkin;Maehara}) also yield monoidal model structures, and the comparison of \cite{Doherty;Kapulkin;Maehara} shows that these models agree with Verity's at the level of bifunctors. Thanks to the work of \cite{Loubaton} (\cite{Gagna;Harpaz;Lanari:Equivalence} in dimension 2) and \cite{Doherty;Kapulkin;Maehara} these models have moreover been compared to models which are universal in the sense of \cite{Barwick;Schommer-Pries:Unicity}. The relevant combinatorics were also studied in detail by \cite{Johnson-Freyd;Scheimbauer} in iterated complete Segal spaces, but that work has not been compared to other models of the Gray tensor product.

More specifically to dimension 2, a bifunctor called the ``Gray tensor product" was constructed for $(\infty,2)$-categories in \cite{Gaitsgory;Rozenblyum} and conjectured to be monoidal biclosed. In \cite{Maehara:Gray} a Gray tensor product on 2-quasicategories is constructed, but the tensor is only lax monoidal at the point-set level, and so it is not immediate that it descends to a monoidal biclosed structure on $\inftytwoCat$. In \cite{Gagna;Harpaz;Lanari:Gray} is constructed a Gray tensor product on scaled simplicial sets which is monoidal biclosed. The underlying bifunctor of this monoidal structure is also shown to coincide at the level of $\infty$-categories with (the 2-dimensional case of) the bifunctor of \cite{Verity:I} (and hence of \cite{Campion;Kapulkin;Maehara}), and steps are taken toward a comparison with the model of \cite{Gaitsgory;Rozenblyum}.

In this work, we leverage the results of \cite{Maehara:Gray} to give a new construction of a Gray tensor product of $(\infty,2)$-categories which comes with a universal property. Namely, we show (\cref{thm:day-refl}) that there is a monoidal biclosed structure on the $\infty$-category of $(\infty,2)$-categories which is uniquely characterized (\cref{cor:univ-prop}) by the fact that on the ``Gray cubes" $\square^n$ (here, $\square^n = [1]^{\otimes n}$ are the Gray tensor powers of the interval), the monoidal structure agrees with the better-understood strict Gray tensor product.
Our proof does not involve very much combinatorics (as long as one is willing to take the results from \cite{Maehara:Gray} for granted) since the only strict $2$-categories we have to deal with are these Gray cubes and the objects in Joyal's category $\Theta_2$, all of which are relatively simple poset-enriched categories.

Unfortunately, our universal property is in some sense ``too strong" in that for a given model of the Gray tensor product, it may be difficult to check that it behaves as expected on the Gray cubes. For example, see \cite[\textsection 7]{Maehara:Gray} for the $2$-quasicategorical case. In particular, we are not able to immediately read off comparison results between our own Gray tensor product and the previously-known models of \cite{Verity:I}, \cite{Campion;Kapulkin;Maehara}, \cite{Gagna;Harpaz;Lanari:Gray}, \cite{Gaitsgory;Rozenblyum}, and \cite{Johnson-Freyd;Scheimbauer}.

%The Gray tensor product we construct was for the most part constructed already in the second author's thesis \cite{Maehara:Gray}, but the role of that work is recast in a new light in the present paper. Here, 
Our construction of the monoidal biclosed structure proceeds as follows.
We begin with the classical Gray tensor product on the category $\square =\{\square^0, \square^1, \square^2, \dots\}$ of Gray cubes, and extend it up to the $\infty$-category $\Psh(\square)$ of presheaves of spaces thereon using the $\infty$-categorical Day convolution (\cref{thm:day-conv}). We then show that $\inftytwoCat$ is canonically a localization of this presheaf category, i.e. that the inclusion $\square \to \inftytwoCat$ is \emph{dense} (\cref{cor:dense-nerve}). In fact we prove something stronger (cf. the first author's \cite{Campion:cubes}): that Joyal's category $\Theta_2$ is contained in the idempotent-completion of $\square$ (\cref{retract}); because $\Theta_2 \to \inftytwoCat$ is well-known to be dense, this implies density of $\square \to \inftytwoCat$ thanks to basic properties of density reviewed in \cref{subsec:infty-categories} below. We do not give explicit generators for this localization. Nevertheless, the work of \cite{Maehara:Gray} allows us to deduce via a formal argument that $\inftytwoCat \subset \Psh(\square)$ is an exponential ideal, so that by the $\infty$-categorical Day reflection theorem, the Day convolution Gray tensor product descends to a monoidal biclosed structure on $\inftytwoCat$ (\cref{thm:day-refl}). The relevant universal properties of Day convolution and Day reflection allow us to read off our model-independent universal property for this Gray tensor product (\cref{cor:univ-prop}).

The outline of the paper is as follows. We begin in \cref{sec:background} with some background. In \cref{subsec:notation} we fix some notation. In \cref{subsec:infty-categories} we recall some facts about density and idempotent completion in the $\infty$-categorical setting. In \cref{subsec:gray-2} we recall the Gray tensor product of strict 2-categories of \cite{Gray}. In \cref{subsec:gray-theta-2} we recall the Gray tensor product of 2-quasicategories from the second author's thesis \cite{Maehara:Gray}. In \cref{sec:dense}, we show that the inclusion $\square \to \inftytwoCat$ is dense (\cref{cor:dense-nerve}) and more specifically that Joyal's category $\Theta_2$ is contained in the idempotent completion of $\square$ (\cref{retract}). This in fact follows (\cref{rmk:campion}) from the main result of \cite{Campion:cubes}, but we prefer to give a direct proof in dimension 2. Finally, in \cref{sec:gray} we construct the Gray tensor product by Day convolution and reflection (\cref{thm:day-conv} and \cref{thm:day-refl}), and read off the desired universal property (\cref{cor:univ-prop}).

\subsection{Acknowledgements}
This project could not have started without Alexander Campbell's initial insight and inspiration which he generously shared with YM. We would like to thank Aaron Mazel-Gee and Dominic Verity for helpful conversations. TC thanks the ARO for their support under MURI Grant W911NF-20-1-0082. YM gratefully acknowledges the support of an Australian Mathematical Society Lift-Off Fellowship and JSPS KAKENHI Grant Number JP21K20329.

\section{Background}\label{sec:background}
This section collects several pieces of background which will be necessary later on. We begin in \cref{subsec:notation} by fixing some notation and conventions, mostly concerning $\infty$-categories and monoidal $\infty$-categories. In the short \cref{subsec:infty-categories}, we recall some facts about density and idempotent completion in the $\infty$-categorical setting. In \cref{subsec:gray-2} we recall the Gray tensor product of strict 2-categories of \cite{Gray}. In \cref{subsec:gray-theta-2} we recall the Gray tensor product of 2-quasicategories from the second author's thesis \cite{Maehara:Gray}.

\subsection{Notation}\label{subsec:notation}
% The term ``$\infty$-category" means ``$(\infty,1)$-category" in this note. We treat $\infty$-categories model-independently, although most of the machinery we use (from \cite{HTT} and \cite{HA}) has been developed in quasicategories. We let $\inftytwoCat$ denote the $\infty$-category of $(\infty,2)$-categories. We also think of $\inftytwoCat$ model-independently, as defined e.g. in \cite{Barwick;Schommer-Pries:Unicity}, though we shall take various nerves of $\inftytwoCat$ to model it via presheaves, especially as presheaves on Joyal's category $\Theta_2$, and so implicitly one may think of our work as taking place in the model of $\Theta_2$-spaces as in \cite{Rezk:cartesian}. However, the technical heart of this work is to apply results (\cite{Maehara:Gray}) from the model of 2-quasicategories of \cite{Ara:nqcat}.

In this note, the term ``$\infty$-category" means ``$(\infty,1)$-category", or more precisely ``quasicategory''; our main references are \cite{HTT} and \cite{HA}.
We let $\inftytwoCat$ denote the $\infty$-category of $(\infty,2)$-categories as defined e.g. in \cite{Barwick;Schommer-Pries:Unicity}.
% Although the statement of our main theorem is model independent (in the sense that it holds for any one of the equivalent $\infty$-categories modelling $\inftytwoCat$) \problematic{If you're being picky about this, you might worry about how model-independent is the concept of a monoidal $\infty$-category}, its proof makes use of the particular model of 2-quasicategories of \cite{Ara:nqcat} so that we can apply the main results from \cite{Maehara:Gray}.

We write $\Spaces$ for the $\infty$-category of spaces. Let $\calC$ be a small $\infty$-category (e.g. $\calC$ could be a small 1-category). We write $\Psh(\calC)$ for the $\infty$-category of $\Spaces$-valued presheaves on $\calC$. We write $\yo : \calC \to \Psh(\calC)$ for the Yoneda embedding. It is not important to us which quasicategory represents $\Psh(\calC)$ up to isomorphism; only up to equivalence. So for instance, $\Psh(\calC)$ might be modeled using left fibrations over $\calC$. We write $\mathrm{N}_\calC : \inftytwoCat \to \Psh(\calC)$ generically for the nerve functor whenever $\calC \subseteq \inftytwoCat$ is a full subcategory (generally when $\calC$ is dense i.e. $\mathrm{N}_\calC$ is fully faithful).

% \problematic{We write $\Spaces$ for the $\infty$-category of spaces. Let $\calC$ be a small $\infty$-category (e.g. $\calC$ could be a small 1-category). We write $\Psh(\calC)$ for the $\infty$-category of $\Spaces$-valued presheaves on $\calC$. We write $\yo : \calC \to \Psh(\calC)$ for the Yoneda embedding. [Define $\nerveinfty$ later; do we need more than one?]}

A \emph{monoidal $\infty$-category $\calC$} is an $\infty$-operad $\calC^\otimes$ cocartesian over the $E_1$-operad in the sense of \cite{HA}. If $\calC, \calD$ are monoidal $\infty$-categories, then we write $\Fun_{E_1}^\lax(\calC, \calD)$ for the $\infty$-category of lax monoidal functors $F :\calC \to \calD$ and monoidal transformations between them (in the notation of \cite{HA}, this would be $\Alg_{\calC/E_1}(\calD)$). We write $\Fun_{E_1}^\strong(\calC,\calD)$ for for the full subcategory of \emph{strong monoidal functors} $F : \calC \to \calD$ (in the notation of \cite{HA}, this would be $\Fun_{E_1}^\otimes(\calC,\calD)$. We let $\Fun_{E_1}^{\lax,\cocts}(\calC, \calD)$ and $\Fun_{E_1}^{\strong,\cocts}(\calC,\calD)$ denote the full subcategories where the functors are required to be colimit-preserving. 

Let $\calC$ be a monoidal biclosed $\infty$-category, and let $\calD \subseteq \calC$ be a full subcategory. We say that $\calD$ is an \emph{exponential ideal} in $\calC$ if, for every $C \in \calC$ and $D \in \calD$, both the left and right internal homs $[C,D]^l$ and $[C,D]^r$ are in $\calD$.
% If $\calC^\otimes, \calD^\otimes$ are $E_1$-monoidal $\infty$-categories, we let $\Fun^{E_1}(\calC^\otimes, \calD^\otimes)$denote the $\infty$-category of lax monoidal functors $\calC^\otimes \to \calD^\otimes$. We let $\Fun^{E_1, \cocts}(\calC^\otimes, \calD^\otimes) \subseteq \Fun^{E_1}(\calC^\otimes, \calD^\otimes)$ denote the full subcategory of functors which preserve small colimits.

We write $\Hom_\calC : \calC^\op \times \calC \to \Spaces$ for the hom-functor of an $\infty$-category $\calC$. We write $\Nat = \Hom_{\Psh(\calC)}$ for an $\infty$-category $\calC$ and we write $\Fun = \Hom_{\inftytwoCat}$.

\subsection{$\infty$-categories}\label{subsec:infty-categories}

Recall that a functor $F: \calC \to \calD$ of $\infty$-categories is said to be \emph{dense} (or \emph{strongly generating} in the terminology of \cite[Definition 4.4.2]{Lurie:Goodwillie} -- ``dense" is preferred as that's the name for the corresponding concept 1-categorically) if the induced functor $F^\ast \calD \to \Psh(\calC)$ is fully faithful. 

\begin{lem}\label{lem:dense-sub}
Let $\calA \subseteq \calB \subseteq \calC$ be inclusions of full sub-$\infty$-categories. Suppose that $\calA \to \calC$ is dense. Then $\calB \to \calC$ is dense.
\end{lem}
\begin{proof}
% This follows from \cref{lem:presheaves-dense}\dots.
This is \cite[Rmk 4.4.7]{Lurie:Goodwillie}.
\end{proof}

\begin{lem}\label{lem:dense-idem}
Let $\calA \to \calB \to \calC$ be inclusions of full sub-$\infty$-categories. Suppose $\calB$ is contained in the idempotent completion of $\calA$. Then $\calB \to \calC$ is dense if and only if $\calA \to \calC$ is dense.
\end{lem}
\begin{proof}
This follows from the fact that the restriction functor $\Psh(\calB) \to \Psh(\calA)$ is an equivalence (cf. \cite[5.1.4.9]{HTT}\footnote{The reference is to the online version of \cite{HTT}.}).
\end{proof}

\begin{rmk}\label{rmk:idem}
Let $\calA$ be an $\infty$-category which happens to be a 1-category. Then the idempotent completion $\tilde \calA$ of $\calA$ is also a 1-category, because each hom-space in $\tilde \calA$ is a retract of a hom-space in $\calA$ (since each object is a retract of a representable in $\Psh(\calA)$ -- see the proof of \cite[5.1.4.2]{HTT}\footnote{The reference is to the online version of \cite{HTT}.}) and hence essentially discrete. Thus $\tilde \calA$ coincides with the idempotent completion of $\calA$ in the 1-categorical sense, comprising the retracts of representables in $\Fun(\calA^\op, \Set)$ (see again the proof of \cite[5.1.4.2]{HTT}).
\end{rmk}

% \begin{definition}
% Let $\calC$ be a monoidal biclosed $\infty$-category, and let $\calD \subseteq \calC$ be a full subcategory. We say that $\calD$ is an \emph{exponential ideal} in $\calC$ if, for every $C \in \calC$ and $D \in \calD$, both the left and right internal homs $[C,D]^l$ and $[C,D]^r$ are in $\calD$.
% \end{definition}

\subsection{Gray tensor product of $2$-categories}\label{subsec:gray-2}
In this subsection, we recall the (\emph{lax}) \emph{Gray tensor product} $\otimes_\square$ on $\twoCat$ \cite[Theorem I.4.9]{Gray}.
The subscript in $\otimes_\square$ indicates that we are considering the Gray tensor product on $\twoCat$ as opposed to those on $\Thetahat$ or $\inftytwoCat$ introduced later.

\begin{definition}
Let $\calA,\calB$ be $2$-categories.
The \emph{Gray tensor product} $\calA \otimes_\square \calB$ is the $2$-category freely generated by the following data:
\begin{itemize}
    \item a $0$-cell $(a,b)$ for each $0$-cell $a$ in $\calA$ and each $0$-cell $b$ in $\calB$,
    \item a $1$-cell $(f,b) : (a,b) \to (a',b)$ for each $1$-cell $f : a \to a'$ in $\calA$ and each $0$-cell $b$ in $\calB$,
    \item a $1$-cell $(a,g) : (a,b) \to (a,b')$ for each $0$-cell $a$ in $\calA$ and each $1$-cell $g : b \to b'$ in $\calB$,
    \item a $2$-cell $(\alpha,b) : (f,b) \to (f',b)$ for each $2$-cell $\alpha : f \to f'$ in $\calA$ and each $0$-cell $b$ in $\calB$,
    \item a $2$-cell $(a,\beta) : (a,g) \to (a,g')$ for each $0$-cell $a$ in $\calA$ and each $2$-cell $\beta : g \to g'$ in $\calB$, and
    \item a $2$-cell $\begin{tikzpicture}[scale = 2, baseline = 25]
        \node (ab) at (0,1) {$(a,b)$};
        \node (ab') at (1.5,1) {$(a,b')$};
        \node (a'b) at (0,0) {$(a',b)$};
        \node (a'b') at (1.5,0) {$(a',b')$};
        \draw[->] (ab) -- node[left] {$(f,b)$} (a'b);
        \draw[->] (ab') -- node[right] {$(f,b')$} (a'b');
        \draw[->] (ab) -- node[above] {$(a,g)$} (ab');
        \draw[->] (a'b) -- node[below] {$(a',g)$} (a'b');
        \draw[->, double] (ab') -- node [above, sloped] {$\gamma_{f,g}$} (a'b);
    \end{tikzpicture}$ for each $1$-cell $f : a \to a'$ in $\calA$ and each $1$-cell $g : b \to b'$ in $\calB$
\end{itemize}
subject to the following conditions:
\begin{itemize}
    \item the assignation $x \mapsto (x,b)$ is a $2$-functor $\calA \to \calA \otimes_\square \calB$ for each $0$-cell $b$ in $\calB$,
    \item the assignation $y \mapsto (a,y)$ is a $2$-functor $\calB \to \calA \otimes_\square \calB$ for each $0$-cell $a$ in $\calA$,
    \item $\begin{tikzpicture}[scale = 2, baseline = 25]
        \node (ab) at (0,1) {$(a,b)$};
        \node (ab') at (1.5,1) {$(a,b')$};
        \node (a'b) at (0,0) {$(a,b)$};
        \node (a'b') at (1.5,0) {$(a,b')$};
        \draw[->] (ab) -- node[left] {$(1_a,b)$} (a'b);
        \draw[->] (ab') -- node[right] {$(1_a,b')$} (a'b');
        \draw[->] (ab) -- node[above] {$(a,g)$} (ab');
        \draw[->] (a'b) -- node[below] {$(a,g)$} (a'b');
        \draw[->, double] (ab') -- node [above, sloped] {$\gamma_{1_a,g}$} (a'b);
    \end{tikzpicture}
    \quad = \quad
    \begin{tikzpicture}[scale = 2, baseline = 25]
        \node (ab) at (0,1) {$(a,b)$};
        \node (ab') at (1.5,1) {$(a,b')$};
        \node (a'b) at (0,0) {$(a,b)$};
        \node (a'b') at (1.5,0) {$(a,b')$};
        \draw[->] (ab) -- node[left] {$1_{(a,b)}$} (a'b);
        \draw[->] (ab') -- node[right] {$1_{(a,b')}$} (a'b');
        \draw[->] (ab) -- node[above] {$(a,g)$} (ab');
        \draw[->] (a'b) -- node[below] {$(a,g)$} (a'b');
        \draw[->, double] (ab') -- node [above, sloped] {$1_{(a,g)}$} (a'b);
    \end{tikzpicture}$
    for each $0$-cell $a$ in $\calA$ and each $1$-cell $g : b \to b'$ in $\calB$,
    \item $\begin{tikzpicture}[scale = 2, baseline = 25]
        \node (ab) at (0,1) {$(a,b)$};
        \node (ab') at (1.5,1) {$(a,b)$};
        \node (a'b) at (0,0) {$(a',b)$};
        \node (a'b') at (1.5,0) {$(a',b)$};
        \draw[->] (ab) -- node[left] {$(f,b)$} (a'b);
        \draw[->] (ab') -- node[right] {$(f,b)$} (a'b');
        \draw[->] (ab) -- node[above] {$(a,1_b)$} (ab');
        \draw[->] (a'b) -- node[below] {$(a',1_b)$} (a'b');
        \draw[->, double] (ab') -- node [above, sloped] {$\gamma_{f,1_g}$} (a'b);
    \end{tikzpicture}
    \quad = \quad
    \begin{tikzpicture}[scale = 2, baseline = 25]
        \node (ab) at (0,1) {$(a,b)$};
        \node (ab') at (1.5,1) {$(a,b)$};
        \node (a'b) at (0,0) {$(a',b)$};
        \node (a'b') at (1.5,0) {$(a',b)$};
        \draw[->] (ab) -- node[left] {$(f,b)$} (a'b);
        \draw[->] (ab') -- node[right] {$(f,b)$} (a'b');
        \draw[->] (ab) -- node[above] {$1_{(a,b)}$} (ab');
        \draw[->] (a'b) -- node[below] {$1_{(a',b)}$} (a'b');
        \draw[->, double] (ab') -- node [above, sloped] {$1_{(f,b)}$} (a'b);
    \end{tikzpicture}$ for each $1$-cell $f : a \to a'$ in $\calA$ and each $0$-cell $b$ in $\calB$,
    \item $\begin{tikzpicture}[scale = 2, baseline = -2]
        \node (ab) at (0,1) {$(a,b)$};
        \node (ab') at (1.5,1) {$(a,b')$};
        \node (a'b) at (0,0) {$(a',b)$};
        \node (a'b') at (1.5,0) {$(a',b')$};
        \node (a''b) at (0,-1) {$(a'',b)$};
        \node (a''b') at (1.5,-1) {$(a'',b')$};
        \draw[->] (ab) -- node[left] {$(f,b)$} (a'b);
        \draw[->] (ab') -- node[right] {$(f,b')$} (a'b');
        \draw[->] (a'b) -- node[left] {$(f',b)$} (a''b);
        \draw[->] (a'b') -- node[right] {$(f',b')$} (a''b');
        \draw[->] (ab) -- node[above] {$(a,g)$} (ab');
        \draw[->] (a'b) -- node[above] {$(a',g)$} (a'b');
        \draw[->] (a''b) -- node[below] {$(a'',g)$} (a''b');
        \draw[->, double] (ab') -- node [above, sloped] {$\gamma_{f,g}$} (a'b);
        \draw[->, double] (a'b') -- node [above, sloped] {$\gamma_{f',g}$} (a''b);
    \end{tikzpicture}
    \quad = \quad 
    \begin{tikzpicture}[scale = 2, baseline = 25]
        \node (ab) at (0,1) {$(a,b)$};
        \node (ab') at (1.5,1) {$(a,b')$};
        \node (a'b) at (0,0) {$(a'',b)$};
        \node (a'b') at (1.5,0) {$(a'',b')$};
        \draw[->] (ab) -- node[left] {$(f'f,b)$} (a'b);
        \draw[->] (ab') -- node[right] {$(f'f,b')$} (a'b');
        \draw[->] (ab) -- node[above] {$(a,g)$} (ab');
        \draw[->] (a'b) -- node[below] {$(a',g)$} (a'b');
        \draw[->, double] (ab') -- node [above, sloped] {$\gamma_{f'f,g}$} (a'b);
    \end{tikzpicture}$ for each composable pair of $1$-cells $f,f'$ in $\calA$ and each $1$-cell $g$ in $\calB$,
    \item $\begin{tikzpicture}[scale = 2, baseline = 25]
        \node (ab) at (0,1) {$(a,b)$};
        \node (ab') at (1.5,1) {$(a,b')$};
        \node (ab'') at (3,1) {$(a,b'')$};
        \node (a'b) at (0,0) {$(a',b)$};
        \node (a'b') at (1.5,0) {$(a',b')$};
        \node (a'b'') at (3,0) {$(a',b'')$};
        \draw[->] (ab) -- node[left] {$(f,b)$} (a'b);
        \draw[->] (ab') -- node[left] {$(f,b')$} (a'b');
        \draw[->] (ab'') -- node[right] {$(f,b'')$} (a'b'');
        \draw[->] (ab) -- node[above] {$(a,g)$} (ab');
        \draw[->] (a'b) -- node[below] {$(a',g)$} (a'b');
        \draw[->] (ab') -- node[above] {$(a,g')$} (ab'');
        \draw[->] (a'b') -- node[below] {$(a',g')$} (a'b'');
        \draw[->, double] (ab') -- node [above, sloped] {$\gamma_{f,g}$} (a'b);
        \draw[->, double] (ab'') -- node [above, sloped] {$\gamma_{f,g'}$} (a'b');
    \end{tikzpicture}
    \quad = \quad
    \begin{tikzpicture}[scale = 2, baseline = 25]
        \node (ab) at (0,1) {$(a,b)$};
        \node (ab') at (1.5,1) {$(a,b'')$};
        \node (a'b) at (0,0) {$(a',b)$};
        \node (a'b') at (1.5,0) {$(a',b'')$};
        \draw[->] (ab) -- node[left] {$(f,b)$} (a'b);
        \draw[->] (ab') -- node[right] {$(f,b')$} (a'b');
        \draw[->] (ab) -- node[above] {$(a,g'g)$} (ab');
        \draw[->] (a'b) -- node[below] {$(a',g'g)$} (a'b');
        \draw[->, double] (ab') -- node [above, sloped] {$\gamma_{f,g'g}$} (a'b);
    \end{tikzpicture}$ for each $1$-cell $f$ in $\calA$ and each composable pair of $1$-cells $g,g'$ in $\calB$,
    \item $\begin{tikzpicture}[scale = 2, baseline = 25]
        \node (ab) at (0,1) {$(a,b)$};
        \node (ab') at (1.5,1) {$(a,b')$};
        \node (a'b) at (0,0) {$(a',b)$};
        \node (a'b') at (1.5,0) {$(a',b')$};
        \draw[->] (ab) .. controls +(-0.4,-0.4) and +(-0.4,0.4) .. node[left, midway] {$(f',b)$} (a'b);
        \draw[->] (ab) .. controls +(0.4,-0.4) and +(0.4,0.4) .. node[right, very near start] {$(f,b)$} (a'b);
        \draw[->] (ab') .. controls +(0.4,-0.4) and +(0.4,0.4) .. node[right, midway] {$(f,b')$} (a'b');
        \draw[->] (ab) -- node[above] {$(a,g)$} (ab');
        \draw[->] (a'b) -- node[below] {$(a',g)$} (a'b');
        \draw[->, double] (ab') -- node [below, sloped] {$\gamma_{f,g}$} (a'b);
        \draw[->, double] (0.2,0.5) -- node [above] {$(\alpha,b)$} (-0.2,0.5);
    \end{tikzpicture}
    \quad = \quad
    \begin{tikzpicture}[scale = 2, baseline = 25]
        \node (ab) at (0,1) {$(a,b)$};
        \node (ab') at (1.5,1) {$(a,b')$};
        \node (a'b) at (0,0) {$(a',b)$};
        \node (a'b') at (1.5,0) {$(a',b')$};
        \draw[->] (ab) .. controls +(-0.4,-0.4) and +(-0.4,0.4) .. node[left, midway] {$(f',b)$} (a'b);
        \draw[->] (ab') .. controls +(0.4,-0.4) and +(0.4,0.4) .. node[right, midway] {$(f,b')$} (a'b');
        \draw[->] (ab') .. controls +(-0.4,-0.4) and +(-0.4,0.4) .. node[left, very near end] {$(f',b')$} (a'b');
        \draw[->] (ab) -- node[above] {$(a,g)$} (ab');
        \draw[->] (a'b) -- node[below] {$(a',g)$} (a'b');
        \draw[->, double] (ab') -- node [above, sloped] {$\gamma_{f,g}$} (a'b);
        \draw[->, double] (1.7,0.5) -- node [above] {$(\alpha,b')$} (1.3,0.5);
    \end{tikzpicture}$ for each $2$-cell $\alpha$ in $\calA$ and each $1$-cell $g$ in $\calB$, and
    \item $\begin{tikzpicture}[scale = 2, baseline = 25]
        \node (ab) at (0,1) {$(a,b)$};
        \node (ab') at (1.5,1) {$(a,b')$};
        \node (a'b) at (0,0) {$(a',b)$};
        \node (a'b') at (1.5,0) {$(a',b')$};
        \draw[->] (ab) -- node[left] {$(f,b)$} (a'b);
        \draw[->] (ab') -- node[right] {$(f,b')$} (a'b');
        \draw[->] (node cs:name=ab, angle=10) .. controls +(0.2,0.2) and +(-0.2,0.2) .. node[midway, above] {$(a,g)$} (node cs:name=ab', angle = 170);
        \draw[->] (node cs:name=a'b, angle=10) .. controls +(0.2,0.2) and +(-0.2,0.2) .. node[very near end, above] {$(a',g)$} (node cs:name=a'b', angle = 170);
        \draw[->] (node cs:name=a'b, angle=-10) .. controls +(0.2,-0.2) and +(-0.2,-0.2) .. node[below] {$(a',g')$} (node cs:name=a'b', angle = -170);
        \draw[->, double] (ab') -- node [above, sloped] {$\gamma_{f,g}$} (a'b);
        \draw[->, double] (0.65,0.1) -- node [right] {$(a',\beta)$} (0.65,-0.1);
    \end{tikzpicture}
    \quad = \quad
    \begin{tikzpicture}[scale = 2, baseline = 25]
        \node (ab) at (0,1) {$(a,b)$};
        \node (ab') at (1.5,1) {$(a,b')$};
        \node (a'b) at (0,0) {$(a',b)$};
        \node (a'b') at (1.5,0) {$(a',b')$};
        \draw[->] (ab) -- node[left] {$(f,b)$} (a'b);
        \draw[->] (ab') -- node[right] {$(f,b')$} (a'b');
        \draw[->] (node cs:name=ab, angle=10) .. controls +(0.2,0.2) and +(-0.2,0.2) .. node[midway, above] {$(a,g)$} (node cs:name=ab', angle = 170);
        \draw[->] (node cs:name=ab, angle=-10) .. controls +(0.2,-0.2) and +(-0.2,-0.2) .. node[very near start, below] {$(a,g')$} (node cs:name=ab', angle = -170);
        \draw[->] (node cs:name=a'b, angle=-10) .. controls +(0.2,-0.2) and +(-0.2,-0.2) .. node[below] {$(a',g')$} (node cs:name=a'b', angle = -170);
        \draw[->, double] (ab') -- node [below, sloped] {$\gamma_{f,g}$} (a'b);
        \draw[->, double] (0.65,1.1) -- node [right] {$(a,\beta)$} (0.65,0.9);
    \end{tikzpicture}$ for each $1$-cell $f$ in $\calA$ and each $2$-cell $\beta$ in $\calB$.
\end{itemize}
\end{definition}

This tensor product is part of a monoidal biclosed structure on $\twoCat$.

\begin{definition}
For $d \ge 0$, the \emph{$d$-dimensional Gray cube} $\square^d$ is the Gray tensor product of $d$ copies of the free-living $1$-cell.
We write $\square$ for the full subcategory of $\twoCat$ spanned by $\square^d$ for $d \ge 0$.
\end{definition}

We can describe $\square^d$ more explicitly as follows.
(See \cite[\textsection 2]{Gray:coherence}, \cite[Theorem 1]{Street:Gray} or \cite[Lemma 3.5]{Maehara:Gray}.)
\begin{itemize}
	\item The objects of $\square^d$ are the binary strings $\eepsilon = (\epsilon_1,\dots,\epsilon_d)\in \{0,1\}^d$ of length $d$.
	\item If two objects $\eepsilon,\zzeta$ satisfy $\epsilon_a \le \zeta_a$ for all $1 \le a \le d$, then the hom-category $\square^d(\eepsilon,\zzeta)$ is the following poset.
	\begin{itemize}
		\item The elements of $\square^d(\eepsilon,\zzeta)$ are precisely the total orders $\pre$ on the set
		\[
		\{a : \epsilon_a = 0, \quad \zeta_a = 1\}.
		\]
		\item A total order $\pre_1$ is ``less than'' $\pre_2$, written $\pre_1 \LHD \pre_2$, if and only if, for any $1 \le a \le b \le d$ with $a \pre_1 b$, we have $a \pre_2 b$.
	\end{itemize}
	Otherwise $\square^d(\eepsilon,\zzeta)$ is empty.
	\item The horizontal composition is given by the join of posets.
	More explicitly, the horizontal composite of $\pre_1$ in $\square^d(\eepsilon,\zzeta)$ and $\pre_2$ in $\square^d(\zzeta,\eeta)$ is the total order $\pre$ given by $a \pre b$ if and only if:
	\begin{itemize}
		\item $a \pre_1 b$;
		\item $a \pre_2 b$; or
		\item $\epsilon_a = 0$, $\zeta_a=\eta_a=1$, $\epsilon_b=\zeta_b = 0$, and $\eta_b = 1$.
	\end{itemize}
\end{itemize}

Intuitively, $\{a : \epsilon_a = 0, \quad \zeta_a = 1\}$ is the set of directions in which one must move to get from $\eepsilon$ to $\zzeta$, and a $1$-morphism $\pre$ simply orders those directions.
A $2$-morphism $\pre_1 \LHD \pre_2$ witnesses that $\pre_2$ is closer  to the usual ordering $\le$ than $\pre_1$ is.

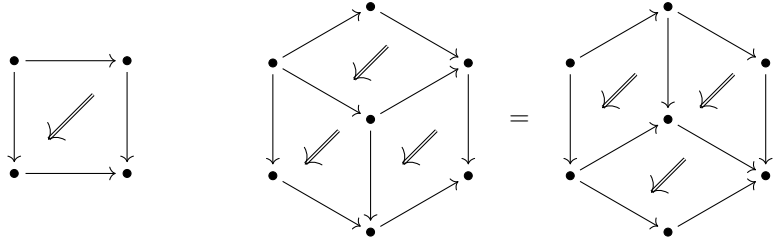
\begin{figure}
\[
\begin{tikzpicture}[scale = 1.5, baseline = 18.5]
	\filldraw
	(0,0) circle [radius = 1pt]
	(1,0) circle [radius = 1pt]
	(0,1) circle [radius = 1pt]
	(1,1) circle [radius = 1pt];
	
	\draw[->] (0.1,0) -- (0.9,0);
	\draw[->] (0.1,1) -- (0.9,1);
	\draw[<-] (0,0.1) -- (0,0.9);
	\draw[<-] (1,0.1) -- (1,0.9);
	
	\draw[->, double] (0.7,0.7) -- (0.3,0.3);
\end{tikzpicture} \hspace{50pt}
\begin{tikzpicture}[baseline = -2,scale = 1.5]	
	\filldraw
	(150:1) circle [radius = 1pt]
	(90:1) circle [radius = 1pt]
	(30:1) circle [radius = 1pt]
	(-30:1) circle [radius = 1pt]
	(-90:1) circle [radius = 1pt]
	(-150:1) circle [radius = 1pt]
	(0,0) circle [radius = 1pt];
	
	\draw[->] (150:1) + (30:0.1) --+ (30:0.9);
	\draw[->] (90:1) + (-30:0.1) --+ (-30:0.9);
	\draw[->] (30:1) + (-90:0.1) --+ (-90:0.9);
	\draw[->] (150:1) + (-90:0.1) --+ (-90:0.9);
	\draw[->] (-150:1) + (-30:0.1) --+ (-30:0.9);
	\draw[->] (-90:1) + (30:0.1) --+ (30:0.9);

	\draw[->] (150:1) + (-30:0.1) --+ (-30:0.9);
	\draw[->] (0:0) + (-90:0.1) --+ (-90:0.9);
	\draw[->] (0:0) + (30:0.1) --+ (30:0.9);
	
	\draw[<-,double] (-150:0.5) + (-0.15,-0.15) --+ (0.15,0.15);
	\draw[<-,double] (90:0.5) + (-0.15,-0.15) --+ (0.15,0.15);
	\draw[<-,double] (-30:0.5) + (-0.15,-0.15) --+ (0.15,0.15);
\end{tikzpicture}
\quad = \quad
\begin{tikzpicture}[baseline = -2,scale = 1.5]
	\filldraw
	(150:1) circle [radius = 1pt]
	(90:1) circle [radius = 1pt]
	(30:1) circle [radius = 1pt]
	(-30:1) circle [radius = 1pt]
	(-90:1) circle [radius = 1pt]
	(-150:1) circle [radius = 1pt]
	(0,0) circle [radius = 1pt];
	
	\draw[->] (150:1) + (30:0.1) --+ (30:0.9);
	\draw[->] (90:1) + (-30:0.1) --+ (-30:0.9);
	\draw[->] (30:1) + (-90:0.1) --+ (-90:0.9);
	\draw[->] (150:1) + (-90:0.1) --+ (-90:0.9);
	\draw[->] (-150:1) + (-30:0.1) --+ (-30:0.9);
	\draw[->] (-90:1) + (30:0.1) --+ (30:0.9);

	\draw[->] (0:0) + (-30:0.1) --+ (-30:0.9);
	\draw[->] (90:1) + (-90:0.1) --+ (-90:0.9);
	\draw[->] (-150:1) + (30:0.1) --+ (30:0.9);
	
	\draw[<-,double] (150:0.5) + (-0.15,-0.15) --+ (0.15,0.15);
	\draw[<-,double] (-90:0.5) + (-0.15,-0.15) --+ (0.15,0.15);
	\draw[<-,double] (30:0.5) + (-0.15,-0.15) --+ (0.15,0.15);
\end{tikzpicture}
\]
\caption{$\square^2$ and $\square^3$}\label{cubes}
\end{figure}

We have depicted $\square^2$ and $\square^3$ in \cref{cubes}; in these pictures, the directions are ordered anticlockwise:
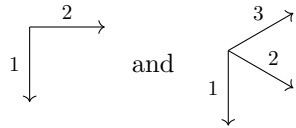
\[
\begin{tikzpicture}[baseline = -17]
	\draw[->] (-3,0) -- (-2,0);
	\draw[->] (-3,0) -- (-3,-1);
	\node[scale = 0.8] at (-2.5,0.2) {$2$};
	\node[scale = 0.8] at (-3.2,-0.5) {$1$};
\end{tikzpicture}
\quad\text{and}\quad
\begin{tikzpicture}[scale = 2, baseline=-8]
	\draw[->] (-2,0) --+ (30:0.5);
	\draw[->] (-2,0) --+ (-30:0.5);
	\draw[->] (-2,0) --+ (-90:0.5);
	\node[scale = 0.8] at (-2.1,-0.25) {$1$};
	\node[scale = 0.8] at (-1.7,-0.05) {$2$};
	\node[scale = 0.8] at (-1.8,0.25) {$3$};
\end{tikzpicture}
\]
respectively.
The $2$-dimensional Gray cube is the ``lax square'' (\cref{cubes} left), and its unique non-trivial $2$-morphism witnesses $(2 \pre 1) \LHD (1 \pre 2)$.
The $3$-dimensional Gray cube is the ``commutative cube'' (\cref{cubes} right), and its long diagonal is the following poset:
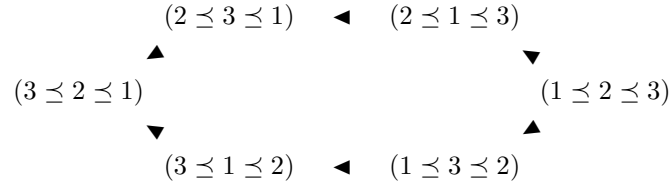
\[
\begin{tikzpicture}
    \node at (-1,0) {$(3 \pre 2 \pre 1)$};
    \node at (6,0) {$(1 \pre 2 \pre 3)$};
    \node at (1,1) {$(2 \pre 3 \pre 1)$};
    \node at (4,1) {$(2 \pre 1 \pre 3)$};
    \node at (1,-1) {$(3 \pre 1 \pre 2)$};
    \node at (4,-1) {$(1 \pre 3 \pre 2)$};
    \node at (2.5,1) {$\LHD$};
    \node at (2.5,-1) {$\LHD$};
    \node[rotate = 30] at (0,0.5) {$\LHD$};
    \node[rotate = -30] at (0,-0.5) {$\LHD$};
    \node[rotate = -30] at (5,0.5) {$\LHD$};
    \node[rotate = 30] at (5,-0.5) {$\LHD$};
\end{tikzpicture}
\]

\subsection{Gray tensor product of $\Theta_2$-sets}\label{subsec:gray-theta-2}

We now recall the Gray tensor product of $\Theta_2$-sets from \cite{Maehara:Gray}.

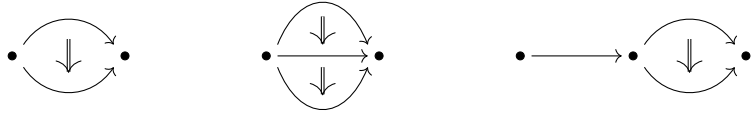
\begin{figure}
\[
\begin{tikzpicture}[baseline = -2, scale = 1.5]
    \filldraw
    (0,0) circle [radius = 1pt]
    (1,0) circle [radius = 1pt];
    \draw [->] (0.1,0.1) .. controls (0.3,0.4) and (0.7,0.4) .. (0.9,0.1);
    \draw [->] (0.1,-0.1) .. controls (0.3,-0.4) and (0.7,-0.4) .. (0.9,-0.1);
    \draw[->, double] (0.5,0.15) -- (0.5,-0.15);
\end{tikzpicture} \hspace{50pt}
\begin{tikzpicture}[baseline = -2, scale = 1.5]
    \filldraw
    (0,0) circle [radius = 1pt]
    (1,0) circle [radius = 1pt];
    \draw [->] (0.1,0.1) .. controls (0.3,0.6) and (0.7,0.6) .. (0.9,0.1);
    \draw [->] (0.1,-0.1) .. controls (0.3,-0.6) and (0.7,-0.6) .. (0.9,-0.1);
    \draw[->] (0.1,0) -- (0.9,0);
    \draw[->, double] (0.5,0.35) -- (0.5,0.1);
    \draw[->, double] (0.5,-0.1) -- (0.5,-0.35);
\end{tikzpicture} \hspace{50pt}
\begin{tikzpicture}[baseline = -2, scale = 1.5]
    \filldraw
    (-1,0) circle [radius = 1pt]
    (0,0) circle [radius = 1pt]
    (1,0) circle [radius = 1pt];
    \draw[->] (-0.9,0) -- (-0.1,0);
    \draw [->] (0.1,0.1) .. controls (0.3,0.4) and (0.7,0.4) .. (0.9,0.1);
    \draw [->] (0.1,-0.1) .. controls (0.3,-0.4) and (0.7,-0.4) .. (0.9,-0.1);
    \draw[->, double] (0.5,0.15) -- (0.5,-0.15);
\end{tikzpicture}
\]
\caption{Objects of $\Theta_2$}\label{Theta}
\end{figure}

\begin{definition}
We write $\nq = [n;q_1,\dots,q_n]$ for the $2$-category freely generated by the $\Cat$-enriched graph such that:
\begin{itemize}
    \item its vertices are $0,1,\dots,n$,
    \item if $0< k \le n$ then $\Hom(k-1,k)$ is the poset $[q_k] = \{0 < 1 < \dots < q_k\}$, and
    \item all other homs are empty.
\end{itemize}
We write $\Theta_2$ for the full subcategory of $\twoCat$ spanned by $\nq$ for $n \ge 0$, $q_1,\dots,q_n \ge 0$.
\end{definition}
For example, (the generating $\Cat$-enriched graphs of) $[1;1]$, $[1;2]$ and $[2;0,1]$ look like \cref{Theta}.

The presheaf category $\Thetahat = [\Theta_2^\op,\Set]$ admits a model structure due to Ara \cite{Ara:nqcat}.
The fibrant objects, called \emph{$2$-quasicategories}, model $(\infty,2)$-categories.
We define a $2$-quasicategorical version of the Gray tensor product as follows.
Write $\nerve : \twoCat \to \Thetahat$ for the nerve functor induced by the inclusion $\Theta_2 \to \twoCat$.

\begin{definition}[{\cite[Definition 3.1]{Maehara:Gray}}]
We define the \emph{$2$-quasicategorical Gray tensor product} $\otimes_{\Theta} : \Thetahat \times \Thetahat \to \Thetahat$ by extending the composite
\[
\Theta \times \Theta \to \twoCat \times \twoCat \xrightarrow{\otimes_\square} \twoCat \xrightarrow{\nerve} \Thetahat
\]
cocontinuously in each variable.
\end{definition}

\begin{thm}[{\cite[Theorem 6.1]{Maehara:Gray}}]\label{thm:quillen}
As a binary functor, $\otimes_\Theta$ is left Quillen with respect to Ara's model structure.
\end{thm}

\begin{cor}\label{cor:two-variable}
The two-variable left adjoint $\otimes_\Theta : \Thetahat \times \Thetahat \to \Thetahat$ induces a two-variable left adjoint $\otimes_\Theta^L : \inftytwoCat \times \inftytwoCat \to \inftytwoCat$.  
\end{cor}
\begin{proof}
By \cite[Theorem 4.6]{Mazel-Gee}, Quillen bifunctors induce two-variable adjunctions of $\infty$-categories, so this follows from \cref{thm:quillen}.
\end{proof}

\begin{rmk}
The functor $\otimes_\Theta^L : \inftytwoCat \times \inftytwoCat \to \inftytwoCat$ will of course be the underlying functor of the monoidal structure we construct on $\inftytwoCat$. However, we will in fact be able to deduce this only indirectly in \cref{rmk:agree} below.
\end{rmk}

\begin{thm}\label{thm:strong}
We have a natural weak equivalence
\[
\begin{tikzpicture}[scale = 2, baseline = 18.5]
        \node (ab) at (0,1) {$\square \times \square$};
        \node (ab') at (1.5,1) {$\Thetahat\times\Thetahat$};
        \node (a'b) at (0,0) {$\square$};
        \node (a'b') at (1.5,0) {$\Thetahat$};
        \draw[->] (ab) -- node[left] {$\otimes_\square$} (a'b);
        \draw[->] (ab') -- node[right] {$\otimes_\Theta$} (a'b');
        \draw[->] (ab) -- node[above] {$\nerve \times \nerve$} (ab');
        \draw[->] (a'b) -- node[below] {$\nerve$} (a'b');
        \draw[->, double] (ab') -- node [above, sloped] {$\sim$} (a'b);
    \end{tikzpicture}
\]
\end{thm}
\begin{proof}
    For $n \ge 0$, write $\otimes^n_\Theta$ for the $n$-ary Gray tensor product on $\Thetahat$ (defined in the same way as $\otimes^2_\Theta = \otimes_\Theta$ but using the $n$-ary Gray tensor product on $\twoCat$ in place of $\otimes_\square$).
    Then for any $X_1,\dots,X_m,Y_1,\dots,Y_n \in \Thetahat$, we have a weak equivalence
    \[
    \otimes^2_\Theta\bigl(\otimes^m_\Theta(X_1,\dots,X_m), \otimes^n_\Theta(Y_1,\dots,Y_n)\bigr) \to \otimes^{m+n}_\Theta(X_1,\dots,X_m,Y_1,\dots,Y_n)
    \]
    by \cite[Corollary 7.11]{Maehara:Gray}.
    In particular, we have a weak equivalence
    \[
    \nerve(\square^m)\otimes_\Theta \nerve(\square^n) \cong \otimes^2_\Theta\bigl(\otimes^m_\Theta(\square^1,\dots,\square^1), \otimes^n_\Theta(\square^1,\dots,\square^1)\bigr) \to \otimes^{m+n}_\Theta(\square^1,\dots,\square^1) \cong \nerve(\square^{m+n}) \cong \nerve(\square^m \otimes_\square \square^n).
    \]
    One can see from the construction of $\otimes_\Theta$ that a $\theta$-cell in $\nerve(\square^m) \otimes_\Theta \nerve(\square^n)$ may be represented (non-uniquely) as a quintuple $(\theta_1,\theta_2,F,F_1,F_2)$ where:
    \begin{itemize}
        \item $\theta_1$ and $\theta_2$ are $2$-categories in $\Theta$,
        \item $F : \theta \to \theta_1 \otimes_\square \theta_2$, $F_1 : \theta_1 \to \square^m$ and $F_2 : \theta_2 \to \square^n$ are $2$-functors.
    \end{itemize}
    The above weak equivalence sends such a $\theta$-cell to the one in $\nerve(\square^m \otimes_\square \square^n)$ corresponding to the composite
    \[
    \begin{tikzcd}
    \theta
    \arrow [r, "F"] &
    \theta_1 \otimes_\square \theta_2
    \arrow [r, "F_1 \otimes_\square F_2"] &
    \square^m \otimes_\square \square^n
    \end{tikzcd}
    \]
    (see \cite[Proposition 7.2]{Maehara:Gray}).
    With this explicit description, its naturality in $\square^m,\square^n \in \square$ is straightforward to verify.
\end{proof}

\section{$\Theta_2$ in Cubes}\label{sec:dense}

In this section, we show (\cref{retract}) that $\Theta_2$ is contained in the idempotent completion of $\square$. The proof unfolds over the course of almost the entire section. At the end, we deduce (\cref{cor:dense-nerve}) that $\square$ is dense in $\inftytwoCat$.

\begin{thm}\label{retract}
Every object of $\Theta_2$ is a retract of an object of $\square$.
\end{thm}

\begin{rmk}
\cref{retract} is a 1-categorical statement. But by \cref{rmk:idem}, it is also an $\infty$-categorical statement.
\end{rmk}

\begin{rmk}\label{rmk:campion}
\cref{retract} may be deduced from a result of the first author \cite[Corollary 2.4]{Campion:cubes}, which says that every object of $\Theta$ is a retract of an object of the $\omega$-dimensional cube category $\square_\omega$. To see this, let $\theta \in \Theta_2$. Then $\theta \in \Theta$. So there is $\square^n \in \square_\omega$ and an idempotent $e : \square^n \to \square^n$ whose splitting is $\theta$. The strict 2-categorical reflection is an idempotent $e_2 : \square^n_2 \to \square^n_2$ on the strict 2-categorical reflection of $\square^n$ whose splitting is the strict 2-categorical reflection of $\theta$. Since $\theta \in \Theta_2$ is already a strict 2-category, this idempotent splitting is $\theta$ itself.

However, we prefer here to give a direct proof in the 2-dimensional case.
\end{rmk}

	Fix $\nq = [n;q_1,\dots,q_n] \in \Theta_2$, and let $d = n+q_1+\dots+q_n$.
	In the following series of lemmas, we prove \cref{retract} by exhibiting $\nq$ as a retract of $\square^d$; that is, we construct $2$-functors $\begin{tikzcd} \nq \arrow [r, "S"] & \square^d \arrow [r, "R"] & \nq \end{tikzcd}$ with $RS = 1$.
	
	We start with the section $S : \nq \to \square^d$.
	\begin{definition}\label{definition:section-rank}
	For any object $k$ in $\nq$, we define its \emph{rank} to be
			\[
			\rho^\Theta(k) = k+q_1+\dots+q_k.
			\]
	\end{definition}
	The object part of $S$ is given by
	\[
	S(k) = (\underbrace{1,\dots,1}_{\rho^\Theta(k)\text{ times}},0,\dots,0).
	\]
	Then for each $i \in [q_k] = \Hom_{\nq}(k-1,k)$, the image $S(i)$ must be a total order $\pre_{k,i}$ on the set of integers $a$ satisfying $\rho^\Theta(k-1)<a\le\rho^\Theta(k)$.
	This total order is obtained as the join of
	\[
	\bigl(\{\rho^\Theta(k-1)+1,\dots,\rho^\Theta(k-1)+i\},\le\bigr)
	\quad\text{and}\quad
	\bigl(\{\rho^\Theta(k-1)+i+1,\dots,\rho^\Theta(k)\}, \ge\bigr).
	\]
	In other words, we have $a \pre_{k,i} b$ if and only if:
	\begin{itemize}
		\item $a \le b \le \rho^\Theta(k-1)+i$;
		\item $a \le \rho^\Theta(k-1)+i < b$; or
		\item $\rho^\Theta(k-1)+i < b \le a$.
	\end{itemize}

\begin{example}
	(a) When $\nq = [1;1]$, we have $\rho^\Theta(0) = 0$ and $\rho^\Theta(1) = 2$.
	The $2$-functor $S$ sends $0 \in [1] = \Hom_{[1;1]}(0,1)$ to $1 \pre 0$ and $1 \in [1]$ to $0 \pre 1$.
	In other words, $S$ simply picks out the interior $2$-cell of the lax square.
	We can visualise the action of $S$ as a distorted version of \cref{Theta} left:
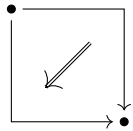
	\[
	\begin{tikzpicture}[scale = 1.5]
		\filldraw
		(1,0) circle [radius = 1pt]
		(0,1) circle [radius = 1pt];
		
		\draw[->] (0,0.9) -- (0,0) -- (0.9,0);
		\draw[->] (0.1,1) -- (1,1) -- (1,0.1);
		
		\draw[->, double] (0.7,0.7) -- (0.3,0.3);
	\end{tikzpicture}
	\]
	
	(b) When $\nq = [1;2]$, we have $\rho^\Theta(0) = 0$, and $\rho^\Theta(1) = 3$.
	The action of $S$ on $\Hom_{[1;2]}(0,1) = [2]$ is given by
	\[
	\begin{split}
	    0 & \mapsto (3 \pre 2 \pre 1)\\
	    1 & \mapsto (1 \pre 3 \pre 2)\\
	    2 & \mapsto (1 \pre 2 \pre 3).
	\end{split}
	\]
	Thus $S$ may be visualised as follows:
	
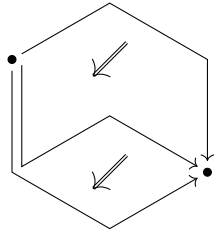
\[
	\begin{tikzpicture}[baseline = -2,scale = 1.5]	
		\filldraw
		(150:1) circle [radius = 1pt]
		(-30:1) circle [radius = 1pt];
		
		\draw[->] (150:1) + (30:0.1) --++ (30:1) --++ (-30:1) --+ (-90:0.9);
		\draw[->] (150:1) + (-90:0.1) --++ (-90:1) --++ (-30:1) --+ (30:0.9);
		\draw[->] (150:1) + (-30:0.1) -- (-150:0.9) -- (0:0) -- (-30:0.9);
		
		\draw[<-,double] (-90:0.5) + (-0.15,-0.15) --+ (0.15,0.15);
		\draw[<-,double] (90:0.5) + (-0.15,-0.15) --+ (0.15,0.15);
	\end{tikzpicture}
\]

(c) When $\nq = [2;0,1]$, $S$ picks out the following cells:
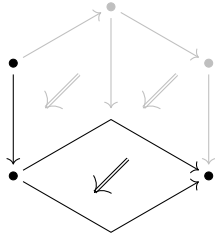
\[
\begin{tikzpicture}[baseline = -2,scale = 1.5]
	\filldraw
	(150:1) circle [radius = 1pt]
	(-30:1) circle [radius = 1pt]
	(-150:1) circle [radius = 1pt];
	
	\filldraw[gray!50!white]
	(90:1) circle [radius = 1pt]
	(30:1) circle [radius = 1pt];
	
	\draw[->,gray!50!white] (150:1) + (30:0.1) --+ (30:0.9);
	\draw[->,gray!50!white] (90:1) + (-30:0.1) --+ (-30:0.9);
	\draw[->,gray!50!white] (30:1) + (-90:0.1) --+ (-90:0.9);
	\draw[->] (150:1) + (-90:0.1) --+ (-90:0.9);
	\draw[->] (-150:1) + (-30:0.1) --++ (-30:1) --+ (30:0.9);
	
	\draw[->,gray!50!white] (90:1) + (-90:0.1) --+ (-90:0.9);
	\draw[->] (-150:1) + (30:0.1) --++ (30:1) --+ (-30:0.9);
	
	\draw[<-,double,gray!50!white] (150:0.5) + (-0.15,-0.15) --+ (0.15,0.15);
	\draw[<-,double] (-90:0.5) + (-0.15,-0.15) --+ (0.15,0.15);
	\draw[<-,double,gray!50!white] (30:0.5) + (-0.15,-0.15) --+ (0.15,0.15);
\end{tikzpicture}
\]
\end{example}
	\begin{lemma}\label{lemma:section}
		The above assignation indeed extends to a unique $2$-functor $S : \nq \to \square^d$.
	\end{lemma}
	\begin{proof}
		It suffices to check that $\pre_{k,i} \LHD \pre_{k,j}$ holds for any $i, j \in [q_k]$ with $i \le j$.
		That is, we must show that $\rho^\Theta(k-1)<a \le b \le \rho^\Theta(k)$ and $a \pre_{k,i} b$ imply $a \pre_{k,j} b$.
		
		Observe that, under the assumption $\rho^\Theta(k-1)<a \le b \le \rho^\Theta(k)$, we have
		\[
		a \pre_{k,i} b \quad \text{iff} \quad a \le \rho^\Theta(k-1)+i
		\]
		and similarly
		\[
		a \pre_{k,j} b \quad \text{iff} \quad a \le \rho^\Theta(k-1)+j.
		\]
		So the desired implication easily follows from the inequality $i \le j$.
	\end{proof}
	
	Next we construct the retraction $R : \square^d \to \nq$.
	\begin{definition}\label{definition:retract-rank}
	For any object $\eepsilon = (\epsilon_1,\dots,\epsilon_d)$ in $\square^d$, we define its \emph{rank} to be
			\[
			\rho^\square(\eepsilon) = \sum_{1 \le a \le d}\epsilon_a.
			\]
	\end{definition}
	We define the object part of $R$ by
	\[
	R(\eepsilon) = \max\bigl\{k : \rho^\Theta(k) \le \rho^\square(\eepsilon)\bigr\}.
	\]
	For each $1 \le a \le d$, let $\ddelta^a \in \square^d$ denote the ``$a$-th unit vector'' given by
	\[
	\delta^a_b = \begin{cases}
		1, & b=a,\\
		0, & \text{otherwise.}
	\end{cases}
	\]
	Then the underlying $1$-category of $\square^d$ is freely generated by the $1$-cells of the form
	\[
	\{a\} : (\eepsilon - \ddelta^a) \to \eepsilon
	\]
	(where $\eepsilon \in \square^d$ satisfies $\epsilon_a = 1$) corresponding to the unique total order on the singleton $\{a\}$.
	Thus, in order to describe the underlying $1$-functor of $R$, it suffices to specify the image of each such $\{a\} : (\eepsilon - \ddelta^a) \to \eepsilon$.
	\begin{itemize}
		\item If $\rho^\square(\eepsilon) = \rho^\Theta(k)$ does not hold for any $k$, then we have $R(\eepsilon - \ddelta^a) = R(\eepsilon)$.
		So the image of $\{a\}$ must be the identity on this object.
		\item Suppose that $\rho^\square(\eepsilon) = \rho^\Theta(k)$ holds for some $k$.
		In this case, the image of $\{a\}$ must be some element of $[q_k]$.
		\begin{itemize}
			\item If $\eepsilon = S(k)$ and moreover $a > \rho^\Theta(k-1)$ then we send $\{a\}$ to $a-\rho^\Theta(k-1)-1 \in [q_k]$.\\
			(Note that $\eepsilon = S(k)$ and $\epsilon_a = 1$ imply $a \le \rho^\Theta(k)$, so $a-\rho^\Theta(k-1)-1$ is indeed a valid element of $[q_k]$.)
			\item Otherwise we send $\{a\}$ to $0 \in [q_k]$.
		\end{itemize}
	\end{itemize}
\begin{example}
	(a) When $\nq = [1;1]$, the $2$-functor $R : \square^2 \to [1;1]$ sends the object $(1,1)$ to $1$ and the other three objects to $0$.
	Of the two atomic $1$-cells with codomain $(1,1)$, it sends the one with domain $(1,0) = (1,1) - \ddelta^2$ to $1$ and that with domain $(0,1) = (1,1) - \ddelta^1$ to $0$.
	We may visualise this $R$ by labelling the $0$- and $1$-cells in \cref{cubes} left with their respective images:
	
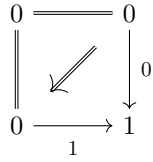
\[
	\begin{tikzpicture}[scale = 1.5]
		\node at (0,0) {$0$};
		\node at (1,0) {$1$};
		\node at (0,1) {$0$};
		\node at (1,1) {$0$};
		
		\draw[->] (0.15,0) -- (0.85,0);
		\draw[double] (0.15,1) -- (0.85,1);
		\draw[double] (0,0.15) -- (0,0.85);
		\draw[<-] (1,0.15) -- (1,0.85);
		
		\draw[->, double] (0.7,0.7) -- (0.3,0.3);
		
		\node[scale = 0.8] at (1.15,0.5) {$0$};
		\node[scale = 0.8] at (0.5,-0.2) {$1$};		
	\end{tikzpicture}
	\]
	(Here an equal sign indicates that $R$ sends the $1$-cell to an identity.)
	
	(b) When $\nq = [1;2]$, $R$ may be depicted as:
	
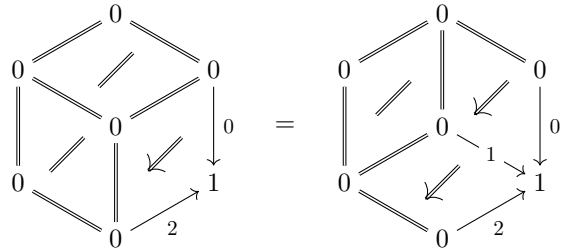
\[
	\begin{tikzpicture}[baseline = -2,scale = 1.5]	
		\node at (150:1) {$0$};
		\node at (90:1) {$0$};
		\node at (-150:1) {$0$};
		\node at (0:0) {$0$};
		\node at (30:1) {$0$};
		\node at (-90:1) {$0$};
		\node at (-30:1) {$1$};
		
		\draw[double] (150:1) + (30:0.15) --+ (30:0.85);
		\draw[double] (90:1) + (-30:0.15) --+ (-30:0.85);
		\draw[->] (30:1) + (-90:0.15) --+ (-90:0.85);
		\draw[double] (150:1) + (-90:0.15) --+ (-90:0.85);
		\draw[double] (-150:1) + (-30:0.15) --+ (-30:0.85);
		\draw[->] (-90:1) + (30:0.15) --+ (30:0.85);

		\draw[double] (150:1) + (-30:0.15) --+ (-30:0.85);
		\draw[double] (0:0) + (-90:0.15) --+ (-90:0.85);
		\draw[double] (0:0) + (30:0.15) --+ (30:0.85);
		
		\draw[double] (-150:0.5) + (-0.15,-0.15) --+ (0.15,0.15);
		\draw[double] (90:0.5) + (-0.15,-0.15) --+ (0.15,0.15);
		\draw[<-,double] (-30:0.5) + (-0.15,-0.15) --+ (0.15,0.15);
		
		\node[scale = 0.8] at (1,0) {$0$};
		\node[scale = 0.8] at (0.5,-0.9) {$2$};
	\end{tikzpicture}
	\quad = \quad
	\begin{tikzpicture}[baseline = -2,scale = 1.5]
		\node at (150:1) {$0$};
		\node at (90:1) {$0$};
		\node at (-150:1) {$0$};
		\node at (0:0) {$0$};
		\node at (30:1) {$0$};
		\node at (-90:1) {$0$};
		\node at (-30:1) {$1$};
		
		\draw[double] (150:1) + (30:0.15) --+ (30:0.85);
		\draw[double] (90:1) + (-30:0.15) --+ (-30:0.85);
		\draw[->] (30:1) + (-90:0.15) --+ (-90:0.85);
		\draw[double] (150:1) + (-90:0.15) --+ (-90:0.85);
		\draw[double] (-150:1) + (-30:0.15) --+ (-30:0.85);
		\draw[->] (-90:1) + (30:0.15) --+ (30:0.85);

		\draw[->] (0:0) + (-30:0.15) --+ (-30:0.85);
		\draw[double] (90:1) + (-90:0.15) --+ (-90:0.85);
		\draw[double] (-150:1) + (30:0.15) --+ (30:0.85);
		
		\draw[double] (150:0.5) + (-0.15,-0.15) --+ (0.15,0.15);
		\draw[<-,double] (-90:0.5) + (-0.15,-0.15) --+ (0.15,0.15);
		\draw[<-,double] (30:0.5) + (-0.15,-0.15) --+ (0.15,0.15);

		\node[scale = 0.8] at (1,0) {$0$};
		\node[scale = 0.8] at (0.5,-0.9) {$2$};
		\node[scale = 0.8, fill = white] at (-30:0.5) {$1$};
	\end{tikzpicture}
	\]
	
	(c) When $\nq = [2;0,1]$, $R$ may be depicted as:
	
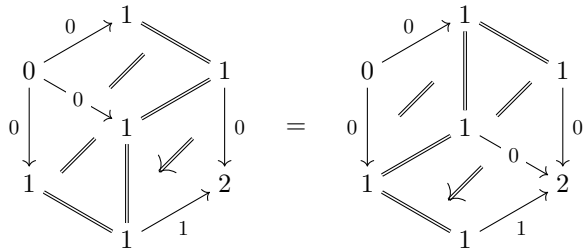
\[
	\begin{tikzpicture}[baseline = -2,scale = 1.5]	
		\node at (150:1) {$0$};
		\node at (90:1) {$1$};
		\node at (-150:1) {$1$};
		\node at (0:0) {$1$};
		\node at (30:1) {$1$};
		\node at (-90:1) {$1$};
		\node at (-30:1) {$2$};
		
		\draw[->] (150:1) + (30:0.15) --+ (30:0.85);
		\draw[double] (90:1) + (-30:0.15) --+ (-30:0.85);
		\draw[->] (30:1) + (-90:0.15) --+ (-90:0.85);
		\draw[->] (150:1) + (-90:0.15) --+ (-90:0.85);
		\draw[double] (-150:1) + (-30:0.15) --+ (-30:0.85);
		\draw[->] (-90:1) + (30:0.15) --+ (30:0.85);

		\draw[->] (150:1) + (-30:0.15) --+ (-30:0.85);
		\draw[double] (0:0) + (-90:0.15) --+ (-90:0.85);
		\draw[double] (0:0) + (30:0.15) --+ (30:0.85);
		
		\draw[double] (-150:0.5) + (-0.15,-0.15) --+ (0.15,0.15);
		\draw[double] (90:0.5) + (-0.15,-0.15) --+ (0.15,0.15);
		\draw[<-,double] (-30:0.5) + (-0.15,-0.15) --+ (0.15,0.15);
		
		\node[scale = 0.8] at (1,0) {$0$};
		\node[scale = 0.8] at (0.5,-0.9) {$1$};
		\node[scale = 0.8] at (-1,0) {$0$};
		\node[scale = 0.8] at (-0.5,0.9) {$0$};
		\node[scale = 0.8, fill = white] at (150:0.5) {$0$};
	\end{tikzpicture}
	\quad = \quad
	\begin{tikzpicture}[baseline = -2,scale = 1.5]
		\node at (150:1) {$0$};
		\node at (90:1) {$1$};
		\node at (-150:1) {$1$};
		\node at (0:0) {$1$};
		\node at (30:1) {$1$};
		\node at (-90:1) {$1$};
		\node at (-30:1) {$2$};
		
		\draw[->] (150:1) + (30:0.15) --+ (30:0.85);
		\draw[double] (90:1) + (-30:0.15) --+ (-30:0.85);
		\draw[->] (30:1) + (-90:0.15) --+ (-90:0.85);
		\draw[->] (150:1) + (-90:0.15) --+ (-90:0.85);
		\draw[double] (-150:1) + (-30:0.15) --+ (-30:0.85);
		\draw[->] (-90:1) + (30:0.15) --+ (30:0.85);

		\draw[->] (0:0) + (-30:0.15) --+ (-30:0.85);
		\draw[double] (90:1) + (-90:0.15) --+ (-90:0.85);
		\draw[double] (-150:1) + (30:0.15) --+ (30:0.85);
		
		\draw[double] (150:0.5) + (-0.15,-0.15) --+ (0.15,0.15);
		\draw[<-,double] (-90:0.5) + (-0.15,-0.15) --+ (0.15,0.15);
		\draw[double] (30:0.5) + (-0.15,-0.15) --+ (0.15,0.15);
		
		\node[scale = 0.8] at (1,0) {$0$};
		\node[scale = 0.8] at (0.5,-0.9) {$1$};
		\node[scale = 0.8] at (-1,0) {$0$};
		\node[scale = 0.8] at (-0.5,0.9) {$0$};
		\node[scale = 0.8, fill = white] at (-30:0.5) {$0$};
	\end{tikzpicture}
	\]
\end{example}
	\begin{lemma}\label{lemma:retract}
		The above assignation indeed extends to a unique $2$-functor $R : \square^d \to \nq$.
	\end{lemma}
	\begin{proof}
		It suffices to check that there is always a (necessarily unique) $2$-cell
		\[
		\begin{tikzpicture}
			\node at (0,2) {$R(\eepsilon-\ddelta^a-\ddelta^b)$};
			\node at (0,0) {$R(\eepsilon-\ddelta^b)$};
			\node at (3.6,2) {$R(\eepsilon-\ddelta^a)$};
			\node at (3.6,0) {$R(\eepsilon)$};
			\draw[->] (0,1.7) -- (0,0.3);
			\draw[->] (3.6,1.7) -- (3.6,0.3);
			\draw[->] (0.9,0) -- (3.1,0);
			\draw[->] (1.35,2) -- (2.7,2);
			\draw[->, double] (2.52,1.4) -- (1.08,0.6);
			\node[scale = 0.8] at (2,-0.3) {$R(\{b\})$};
			\node[scale = 0.8] at (2,2.3) {$R(\{b\})$};
			\node[scale = 0.8] at (-0.6,1) {$R(\{a\})$};
			\node[scale = 0.8] at (4.2,1) {$R(\{a\})$};
		\end{tikzpicture}
		\]
	for any $\eepsilon \in \square^d$ and $1 \le a < b \le d$ with $\epsilon_a = \epsilon_b = 1$.
	Note that, since the object part of $R$ is solely determined by the rank, we have $R(\eepsilon-\ddelta^a) = R(\eepsilon-\ddelta^b)$.
	Thus we may instead exhibit $2$-cells of the form
	\[
	\begin{tikzpicture}
		\node at (0,2) {$R(\eepsilon-\ddelta^a-\ddelta^b)$};
		\node at (0,0) {$R(\eepsilon-\ddelta^b)$};
		\node at (3.6,2) {$R(\eepsilon-\ddelta^a)$};
		\node at (3.6,0) {$R(\eepsilon)$};
		\draw[->] (0,1.7) -- (0,0.3);
		\draw[->] (3.6,1.7) -- (3.6,0.3);
		\draw[->] (0.9,0) -- (3.1,0);
		\draw[->] (1.35,2) -- (2.7,2);
		\draw[double] (0.54,0.3) -- (3.06,1.7);
		\node[scale = 0.8] at (2,-0.3) {$R(\{b\})$};
		\node[scale = 0.8] at (2,2.3) {$R(\{b\})$};
		\node[scale = 0.8] at (-0.6,1) {$R(\{a\})$};
		\node[scale = 0.8] at (4.2,1) {$R(\{a\})$};
		\draw[->, double] (1,1.5) -- (1,1);
		\draw[->, double] (2.6,1) -- (2.6,0.5);
	\end{tikzpicture}
	\]
	
	First we consider the upper triangle.
	Note that if $\rho^\square(\eepsilon - \ddelta^a) = \rho^\Theta(k)$ does not hold for any $k$, then we trivially have such a $2$-cell since all three edges of the triangle are the identity.
	So assume that we have $\rho^\square(\eepsilon - \ddelta^a) = \rho^\Theta(k)$ for some $k$.
	Observe that $\eepsilon - \ddelta^a$ cannot be $S(k)$ since its $a$-th coordinate is $0$ and its $b$-th coordinate is $1$.
	Thus the top edge $R(\{b\})$ is $0 \in [q_k]$, which implies the existence of the desired $2$-cell.
	
	For the lower triangle, we consider the following cases separately.
	\begin{itemize}
		\item If $\rho^\square(\eepsilon) = \rho^\Theta(k)$ does not hold for any $k$, then we trivially obtain the desired $2$-cell since all three edges of the triangle are the identity.
		\item If $\rho^\square(\eepsilon) = \rho^\Theta(k)$ but $\eepsilon \neq S(k)$ then we trivially obtain the desired $2$-cell since $R(\{a\}) = R(\{b\}) = 0$.
		\item Suppose that $\eepsilon = S(k)$ holds for some $k$.
		\begin{itemize}
			\item If $a \le \rho^\Theta(k-1)$ then $R(\{a\}) = 0 \in [q_k]$ is the initial object and so we obtain the desired $2$-cell.
			\item If $a > \rho^\Theta(k-1)$ then the inequality
			\[
			R(\{a\}) = a-\rho^\Theta(k-1)-1 > b-\rho^\Theta(k-1)-1 = R(\{b\})
			\]
			corresponds to the desired $2$-cell.
		\end{itemize}
	\end{itemize}
	This completes the proof.
\end{proof}

\begin{lemma}\label{lemma:section-retract}
	The composite $RS$ is the identity on $\nq$.
\end{lemma}
\begin{proof}
	That it is the identity on the $0$-cells is easy to check.
	For the $1$-cells, fix $i \in [q_k]$, and consider the decomposition of its image $S(i)$ into atomic $1$-cells in $\square^d$.
	One can easily check (by considering the rank of the objects appearing in the decomposition) that $R$ sends each factor except for the last to the identity at $k-1$.
	This last factor is
	\[
	\{\rho^\Theta(k-1)+i+1\} : \bigl(S(k)-\ddelta^{\rho^\Theta(k-1)+i+1}\bigr) \to S(k),
	\]
	and $R$ sends it back to
	\[
	\bigl(\rho^\Theta(k-1)+i+1\bigr)-\rho^\Theta(k-1)-1 = i.
	\]
	This completes the proof.
\end{proof}

\begin{proof}[Proof of \cref{retract}]
We exhibit $\theta \in \Theta_2$ as a retract of $\square^d \in \square$. The section is given after \cref{definition:section-rank}. The retract is given after \cref{definition:retract-rank}. The correctness of the construction is given by \cref{lemma:section}, \cref{lemma:retract}, and \cref{lemma:section-retract}.
\end{proof}

\begin{cor}\label{cor:dense-nerve}
The fully faithful functor $\square \to \inftytwoCat$ is dense. Explicitly, the restricted Yoneda embedding $\inftytwoCat \to \Psh(\square)$ is fully faithful, with a left adjoint $\refl$.
\end{cor}
\begin{proof}
We have $\Theta_2 \subseteq \tilde \square \supseteq \square$ (where $\tilde \square$ is the idempotent completion of $\square$), inducing restriction functors $\Psh(\Theta_2) \leftarrow \Psh(\tilde \square) \to \Psh(\square)$. As in \cref{lem:dense-idem}, $\Psh(\tilde \square) \to \Psh(\square)$ is an equivalence. By \cref{lem:dense-sub} and \cref{retract}, $\Psh(\Theta_2) \leftarrow \Psh(\tilde \square)$ is an equivalence. It is well-known that $\Theta_2$ is dense in $\inftytwoCat$ (for example, $\inftytwoCat$ may be presented via the $\Theta$-space model), i.e. the nerve functor $\inftytwoCat \to \Psh(\Theta_2)$ is fully faithful. Composing with the equivalence $\Psh(\Theta_2) \simeq \Psh(\tilde \square) \simeq \Psh(\square)$, we find that the nerve functor $\inftytwoCat \to \Psh(\square)$ is fully faithful. The left adjoint $\refl$ is likewise obtained by composition with the localization from $\Psh(\Theta_2)$ to $\inftytwoCat$.
\end{proof}

\section{The Gray Tensor Product}\label{sec:gray}

In this section, we use the results of the previous section to construct a Gray tensor product on $\inftytwoCat$ (\cref{thm:day-refl}). The construction is in two steps. First, we extend $\otimes_\square$ by Day convolution from $\square$ to $\Psh(\square)$ (\cref{thm:day-conv}). We then use the good properties of $\otimes_\Theta$ recalled in \cref{subsec:gray-theta-2} to show that this Day convolution tensor product descends to $\inftytwoCat$, which a full subcategory of $\Psh(\square)$ by \cref{cor:dense-nerve}. Finally, we deduce a universal property of the Gray tensor product (\cref{cor:univ-prop}). Given the results of the previous sections and of \cite{HA}, the arguments in this section are purely formal.

\begin{thm}\label{thm:day-conv}
There is a monoidal biclosed structure $\hat \otimes_\square$ on $\Psh(\square)$ such that the Yoneda embedding $\yo : \square \to \Psh(\square)$ is strong monoidal. Moreover, for any cocomplete monoidal $\infty$-category $\calE$ with colimits preserved separately in each variable, $\yo$ induces an equivalence $\Fun_{E_1}^{\lax,\cocts}(\Psh(\square),\calE) \to \Fun_{E_1}^\lax(\square,\calE)$.
\end{thm}
\begin{proof}
This is \cite[Proposition 4.8.1.10]{HA}, applied in the case where $\calK = \emptyset$, $\calK'$ consists of all small $\infty$-categories, and $\calO$ is the $E_1$ operad, and $\calC = \square$, considered as a monoidal category under $\otimes_\square$. (Strictly speaking, we deduce from here the above statement amended to say that $\hat \otimes_\square$ preserves colimits separately in each variable, and we deduce that $\Psh(\square)$ is monoidal biclosed by the $\infty$-categorical adjoint functor theorem (\cite[Corollary 5.5.2.9]{HTT}).)
\end{proof}

\begin{lem}\label{lem:ideal}
Let $\begin{tikzcd}
L : \calC \ar[r,shift left] & \calD \ar[l, shift left] : R
\end{tikzcd}$
be an adjunction $L \dashv R$ of $\infty$-categories, where $R$ is fully faithful. Let $T : \calC \times \calC \to \calC$ be a bifunctor. Suppose that there is a functor $H : \calC^\op \times \calD \to \calD$ such that there is an equivalence $\Hom_\calC(C', H(C, D)) \cong \Hom_\calC(T(C',C), D)$, natural in $C',C,D$. Then $T$ is left-compatible with the localization $L$ in the sense that if $f : C'' \to C'$ is carried by $L$ to an equivalence for $C'',C' \in \calC$, then $f \otimes \id_{C'''}$ is carried by $L$ to an equivalence for $C''' \in \calC$.
\end{lem}
\begin{proof}
It suffices to check that $\Hom_\calC(f \otimes \id_{C'''}, D) : \Hom_\calC(T(C', C'''), D) \to \Hom_\calC(T(C'',C'''), D)$ is an equivalence for $D \in \calD$. We have 
\[
\Hom_\calC(T(C',C'''), D) \cong \Hom_\calC(C', H(C''', D)) \cong \Hom_\calC(C'', H(C''', D)) \cong \Hom_\calC(T(C'', C'''), D)
\]
as desired.
\end{proof}

\begin{thm}\label{thm:day-refl}
The fully faithful nerve $\inftytwoCat \to \Psh(\square)$ exhibits $\inftytwoCat$ as an exponential ideal in the monoidal category $\Psh(\square)$ of \cref{thm:day-conv}. 
% There is a unique monoidal biclosed structure $\otimes$ on $\inftytwoCat$ such that $\nerveinfty$ is lax monoidal. 
Thus the left adjoint $\refl$ to $\nerveinfty$ is strong monoidal. Moreover, for any cocomplete monoidal $\infty$-category $\calE$ with colimits preserved by $\otimes$ separately in each variable, precomposition with $\refl$ induces a fully faithful functor
\[
\Fun_{E_1}^{\strong,\cocts}(\inftytwoCat, \calE) \to \Fun_{E_1}^{\strong, \cocts}(\Psh(\square), \calE)
\]
whose essential image comprises those strong monoidal, colimit preserving functors $\Psh(\square) \to \calE$ whose underlying functor descends to a functor $\inftytwoCat \to \calE$. 
\end{thm}
\begin{proof}
By 
% \cite[Proposition 2.2.1.9]{HA}
\cite[Proposition 4.1.7.4]{HA}, it suffices to show that the $E_1$ monoidal structure $\hat \otimes_\square$ is compatible with the localization functor $L$. That is, we must show that if $f : X \to Y$ is a map of presheaves such that $L(f)$ is an equivalence, then $L(f \hat \otimes_\square \id_Z)$ and $L(\id_Z \hat \otimes_\square f)$ are equivalences for each presheaf $Z$. We check the first statement; the proof in the other case is similar. Since $L$ commutes with colimits and $\Psh(\square)$ is generated under colimits by representables, it suffices to check this when $Z = \yo z$ is representable. 

By \cref{lem:ideal}, it will suffice to show that there is a functor $H : \Psh(\square)^\op \times \inftytwoCat \to \inftytwoCat$ and a natural equivalence $\Nat(X, H(Y, D)) \cong \Nat(X \hat \otimes_\square Y, D)$ for $X, Y \in \Psh(\square)$ and $D \in \inftytwoCat$. By \cref{cor:two-variable}, there is a functor $[-,-] : \inftytwoCat^\op \times \inftytwoCat \to \inftytwoCat$ with a natural equivalence $\Fun(A \otimes_\Theta^L B, D) \cong \Fun(A, [B, D])$. We restrict this to a functor $\square^\op \times \inftytwoCat \to \inftytwoCat$, and then Kan extend to a functor $H : \Psh(\square)^\op \times \inftytwoCat \to \inftytwoCat$. We have
\begin{align*}
    \Nat(\yo \square^m, H(\yo \square^n, \nerveinfty D))
    &\cong \Nat(\yo \square^m, \nerveinfty [\nerve\square^n, D]) \qquad \text{by definition of $H$}\\
    &\cong \Fun(\nerve\square^m, [\nerve\square^n, D]) \qquad \text{by \cref{cor:dense-nerve}}\\
    &\cong \Fun(\nerve\square^m \otimes_\Theta^L \nerve\square^n, D) \qquad \text{by \cref{cor:two-variable}} \\
    &\cong \Fun(\nerve\square^m \otimes_\Theta \nerve\square^n, D) \qquad \text{by \cref{thm:quillen} (every object of $\Thetahat$ is cofibrant.)} \\
    &\cong \Fun(\nerve(\square^m \otimes_\square \square^n), D) \qquad \text{by \cref{thm:strong}} \\
    &\cong \Nat(\yo (\square^m \otimes_\square \square^n), \nerveinfty D) \qquad \text{by \cref{cor:dense-nerve}}\\
    &\cong \Nat((\yo \square^m) \hat \otimes_\square (\yo \square^n), \nerveinfty D) \qquad \text{by \cref{thm:day-conv}}
\end{align*}
naturally in $\square^m, \square^n \in \square$ and $D \in \inftytwoCat$. As both sides are limit-preserving in all variables, this extends by \cref{cor:dense-nerve} to an equivalence $\Nat(X, H(Y, \nerveinfty D)) \cong \Nat(X \hat \otimes_\square Y, \nerveinfty D)$ natural in $X,Y \in \Psh(\square), D \in \inftytwoCat$.
\end{proof}

\begin{rmk}\label{rmk:agree}
The underlying bifunctor $\otimes$ of the monoidal structure of \cref{thm:day-refl} is none other than the bifunctor $\otimes_\Theta^L: \inftytwoCat \times \inftytwoCat \to \inftytwoCat$ of \cref{cor:two-variable}. To see this, observe that $\otimes_\Theta^L$ and $\hat \otimes_\square$ both commute with colimits in each variable, and agree with $\otimes_\square$ on cubes, which are dense in $\inftytwoCat$ by \cref{cor:dense-nerve}.
\end{rmk}

\begin{cor}\label{cor:univ-prop}
The fully faithful functor $\square \to \inftytwoCat$ is strong monoidal. Here $\square$ is regarded as a monoidal ($\infty$-)category under the usual Gray tensor product $\otimes_\square$ of strict 2-categories, and $\inftytwoCat$ is regarded as monoidal under the Gray tensor product $\otimes$ of \cref{thm:day-refl}. Moreover, we have the following two universal properties:
\begin{enumerate}
    \item For any monoidal biclosed, cocomplete $\infty$-category $\calE$, composition with $\square \to \inftytwoCat$ induces a fully faithful functor 
    \[\Fun_{E_1}^{\strong, \cocts}(\inftytwoCat, \calE) \to \Fun_{E_1}^\strong(\square, \calE)\]
    whose essential image comprises those monoidal functors $\square \to \calE$ whose underlying functor extends to a cocontinuous functor $\inftytwoCat \to \calE$.
    \item $\otimes$ is the unique monoidal biclosed structure on $\inftytwoCat$ such that $\square \to \inftytwoCat$ is strong monoidal.
\end{enumerate}
\end{cor}
\begin{proof}
We can factorise the functor as $\square \xrightarrow{\yo} \Psh(\square) \xrightarrow{\refl} \inftytwoCat$.
The first factor is strong monoidal by \cref{thm:day-conv} and the second factor is strong monoidal by \cref{thm:day-refl}. The first universal property follows by concatenating \cref{thm:day-conv} and \cref{thm:day-refl}. For the second, suppose that $F : \square \to \inftytwoCat$ is strong monoidal with respect to another monoidal biclosed structure $\otimes'$ on $\inftytwoCat$, and that its underlying functor is the usual $\square \to \inftytwoCat$. By the first universal property, there is an induced cocontinuous, monoidal biclosed functor $\hat F : \inftytwoCat \to \inftytwoCat$ which is strong monoidal from $\otimes$ to $\otimes'$, and whose restriction to $\square$ is the identity. By \cref{cor:dense-nerve}, $\hat F$ is naturally isomorphic to the identity. So $\hat F$ is an equivalence. Since $\hat F$ is strong monoidal, it is a monoidal equivalence \cite[Rmk 2.1.3.8]{HA}.
\end{proof}

\bibliographystyle{alpha}
\bibliography{gray}

\end{document}